\documentclass[final,onefignum,onetabnum]{siamart250211}

\usepackage{lipsum}
\usepackage{amsfonts}
\usepackage{graphicx}
\usepackage{epstopdf}
\usepackage{algorithmic}
\ifpdf
  \DeclareGraphicsExtensions{.eps,.pdf,.png,.jpg}
\else
  \DeclareGraphicsExtensions{.eps}
\fi
\usepackage{cprotect}
\usepackage{subcaption}

\newsiamremark{remark}{Remark}
\newsiamremark{hypothesis}{Hypothesis}
\crefname{hypothesis}{Hypothesis}{Hypotheses}
\newsiamthm{claim}{Claim}
\newsiamremark{fact}{Fact}
\crefname{fact}{Fact}{Facts}

\headers{Computing the zeros of a holomorphic function}{Bowhay, Nakatsukasa, and Zaid}

\title{Computing the zeros of a holomorphic function using quadrature-based subdivision and rational approximation of the logarithmic derivative}

\author{Jake Bowhay\thanks{\email{jake.bowhay@bristol.ac.uk}, School of Engineering Mathematics and Technology, University of Bristol, Bristol, BS8 1TW \& Mathematical Institute, University of Oxford, Oxford, OX2 6GG} \and
Yuji Nakatsukasa \thanks{\email{nakatsukasa@maths.ox.ac.uk}, Mathematical Institute, University of Oxford, Oxford, OX2 6GG}\and
Irwin Zaid \thanks{\email{irwin.zaid@chch.ox.ac.uk}, Christ Church,
University of Oxford, Oxford, OX1 1DP}
}

\usepackage{amsopn}

\DeclareMathOperator{\res}{Res}
\DeclareMathOperator{\tr}{tr}
\DeclareMathOperator{\adj}{adj}

\ifpdf
\hypersetup{
  pdftitle={Computing the zeros of a holomorphic function using quadrature-based subdivision, the AAA algorithm, and the logarithmic derivative.},
  pdfauthor={J. Bowhay, Y. Nakatsukasa, I. Zaid}
}
\fi

\newcommand*\diff{\mathop{}\!\mathrm{d}}

\begin{document}

\maketitle

% REQUIRED
\begin{abstract}
 We introduce a new method that uses AAA approximation to reliably compute all the zeros of a holomorphic function in a specified search region in the complex plane. Specifically, the method is based on rational approximation of the logarithmic derivative in combination with subdivision of the search region based on Cauchy's argument principle. This is motivated by the fact that, while it is straightforward to compute the zeros of a AAA rational approximation, there is no guarantee that all of the zeros of the function being approximated will be found. Many of the ideas presented are also applicable to computing both the zeros and the poles of a meromorphic function. An implementation of the method is provided by the Python package \verb|skzeros|.
\end{abstract}

% REQUIRED

% REQUIRED

\section{Introduction}
In this paper, we present a method for solving the following problem: given a function $f:W\rightarrow\mathbb{C}$, where $f(z)$ is holomorphic in $W\subseteq\mathbb{C}$, the function's derivative $f^\prime(z)$ (or a numerical approximation), and a simply connected, bounded search region $\Omega\subset W$, we wish to compute \emph{all} of the zeros of $f(z)$ that lie within $\Omega$. This method can also find \emph{all} the poles of a meromorphic function $f(z)$ provided that $1/f(z)$ is holomorphic in $\Omega$. If this is not the case (e.g. if $f$ is meromorphic with both poles and zeros) many of the ideas from the method are still applicable, however, it is no longer possible to confirm that all of the zeros and poles in the search region have been found.

Computing the zeros and poles of a function is a fundamental task in solving many mathematical problems. For example, zero finding occurs in applications such as simulating the collision of solitons~\cite{sinegordon}, the analysis of waves in an antenna~\cite{Long1997}, and the study of multiplicity-induced-dominancy for scalar delay-differential equations~\cite{dde}. Likewise, pole finding has widespread applications, from inverting the spectral Green's function of the Helmholtz problem~\cite{Olyslager2004} to solving problems in fractional calculus~\cite{Sheng2011}. Locating the poles of a function also has applications beyond the residue theorem. For example, it is possible to reformulate generalized eigenvalue~\cite{Sakurai2003} and nonlinear eigenvalue~\cite{Asakura2009,AAAsubdivision,Gttel2022} problems into pole-finding problems.

Traditionally, numerical computation of the zeros of a function is done using an iterative root-finding method, such as Newton's or Muller's method~\cite{Muller1956}. These methods require the user to supply an initial starting guess and then iteratively converge to a single zero of the function. Therefore, one approach to finding all the zeros in a region is to use an iterative method with a range of initial guesses. This can be an effective method if approximate locations of the zeros are known \emph{a priori} as they can be used as an initial guess.

When the approximate locations are unknown, this approach has several difficulties. If the number of zeros in the search region is unknown, it is difficult to know when to stop searching for further zeros. Furthermore, if the initial guess is near a stationary point, the method may take a very large step, converging to a zero outside the region of interest. In fact, predicting which zero Newton's method will converge to is not a straightforward task, as illustrated by the so-called Newton fractal~\cite{Epureanu1998}.

Another issue is that it is likely that multiple initial guesses will converge to the same zero, which is problematic as it greatly increases the number of attempts needed to find all of the zeros. One approach to overcome this is numerical deflation, as described in Chapter 2 of~\cite{Wilkinson1994-cd} or~\cite{Farrell2015}. Theoretically, this makes it much easier to find all the zeros in a region, as performing deflation after each zero is found means that the iterative method will not converge to that zero again. Extensions of this technique have been applied to various problems, see~\cite{Brow1971,Farrell2015}. Nevertheless, deflation does not guarantee that all solutions will be found~\cite{Farrell2015}. Furthermore, the error in the computed zeros means they are not exactly factored out when we perform deflation. This means that deflation may not always work well in practice, especially if there are repeated zeros, see Chapter~4.4 of~\cite{Johnson1977-aw}.

As identified by Delves and Lyness~\cite{Delves1967}, the limitations of local iterative methods motivate an alternative approach in which all the zeros of a region are found at once. For example, when the search region is restricted to an interval then \verb|Chebfun|~\cite{Battles2004} can quickly and accurately compute all the zeros of a function. In this paper, we present a new method of this type, which first subdivides the search region based on the value of Cauchy's argument principle and then approximates the logarithmic derivative $f^\prime(z)/f(z)$ in each of the subdivided regions using the adaptive Antoulas–Anderson (AAA) algorithm~\cite{Nakatsukasa2018} in order to compute the location of the zeros.

\subsection{Outline}

The structure of the paper is as follows. In \Cref{sec:quad_methods}, we discuss existing quadrature-based methods as they introduce the idea of quadrature-based subdivision of the search region. \Cref{sec:AAA} briefly introduces the AAA algorithm. \Cref{sec:logderiv} shows how approximating the logarithmic derivative opposed to the function $f(z)$ itself can improve the performance of AAA for zero-finding and discusses how to approximate the derivative if it is not explicitly known. \Cref{sec:subdiv} shows the benefits of using a divide-and-conquer approach in which the search region is subdivided into smaller subregions before applying AAA. \Cref{sec:brach} shows how to handle functions that have branch cuts within the search region. \Cref{sec:comp} compares the AAA algorithm to existing approaches and shows that AAA requires orders of magnitudes fewer function evaluations. \Cref{sec:AAA_Quad} presents a method that combines quadrature-based subdivision of the search region with AAA approximation of the logarithmic derivative. Finally, \Cref{sec:Num_expr} presents some numerical experiments demonstrating the method and \Cref{sec:conclusion} gives some concluding remarks.

\subsection{Code Availability}
\label{sec:code}
A Python implementation of the method described in \Cref{sec:AAA_Quad} and code to reproduce all of the numerical experiments and figures presented in this paper are available at \url{https://github.com/j-bowhay/skzeros} or as part of the Python package \verb|skzeros|.

\section{Quadrature-Based Methods}
\label{sec:quad_methods}

In this section, we give an overview of two existing methods for computing all of the zeros of a holomorphic function in a given region that introduce ideas relevant to the method presented later in this paper. These are the method introduced by Delves and Lyness \cite{Delves1967} and its subsequent improvement by Kravanja and Van Barel~\cite{Kravanja2000}. For a more complete overview of other zero-finding methods see~\cite{keyMcNameea,keyMcNameeb}. The classic 1967 paper by Delves and Lyness~\cite{Delves1967} introduces a method to locate all the zeros of a given holomorphic function $f(z)$ that lie in the interior of a region $\Omega\subset\mathbb{C}$ using numerical quadrature. Here, we give a summary of the method.

The method is based on the following observation: if the positively oriented boundary $\partial\Omega$ of the search region $\Omega$ does not pass through any zeros of $f(z)$ then we have the following result for the $k$th ordinary moment of the argument principle
\begin{equation}
    s_k = \frac{1}{2\pi i}\oint_{\partial\Omega} z^k \frac{f^\prime(z)}{f(z)} \diff{z} = \sum_{j=1}^n \alpha_ja_j^k,
    \label{eq:delves}
\end{equation}
where $i^2 = -1$ and $a_1,\dots,a_n$ are the zeros of $f(z)$ which lie in $\Omega$ with associated multiplicities $\alpha_1,\dots,\alpha_n$. This result follows directly from the generalised argument principle for a meromorphic function $f(z)$
\begin{equation}
        \frac{1}{2\pi i}\oint_{\partial\Omega}g(z)\frac{f^\prime(z)}{f(z)}  \diff{z} = \sum_{j=1}^n \alpha_j g(a_j) - \sum_{j=1}^m \beta_j g(b_j),
\end{equation}
where $b_1,\dots,b_m$ are the poles of $f(z)$ which lie in $\Omega$ with associated multiplicities $\beta_1,\dots,\beta_m$. Here, $g(z)\equiv z^k$ gives \eqref{eq:delves} and $g(z)\equiv1$ recovers the familiar Cauchy's argument principle.

By computing the result of the contour integral \eqref{eq:delves} using numerical quadrature, the resulting values of $s_0,s_1,s_2,\dots,s_N$, where $N=s_0$ is the number of zeros in $\Omega$, can be used to generate a polynomial $p(z)$ of degree $N$. This polynomial is called the equivalent polynomial and has zeros at the same locations as $f(z)$ in $\Omega$ with matching multiplicities. The computation of the equivalent polynomial is described later by \eqref{eq:equiv_poly}. We note that the equivalent polynomial is \emph{not} a polynomial approximation to $f(z)$ despite sharing the same zeros in $\Omega$. If $f(z)$ has many zeros in $\Omega$ then the polynomial $p(z)$ is of high degree and both the map from ordinary moments $s_k$ to coefficients of the equivalent polynomial and the map from coefficients to zeros of $p(z)$ have the potential to be highly ill-conditioned, see Chapter 1 of~\cite{Kravanja2000}. Delves and Lyness' key idea is to subdivide $\Omega$ into smaller regions with a manageable number of zeros to avoid this ill-conditioning. This method can be summarized by the following steps:
\begin{enumerate}
    \item Determine the number of zeros $N=s_0$ in the region. If $N$ is less than a user-supplied parameter (typically 5), then $s_1,s_2,\dots,s_N$ are computed, and the method proceeds to Step~3. If not the method proceeds to Step~2.
    \item Subdivide the region into smaller regions and perform Step~1 on each of the new subregions.
    \item Once all regions contain few enough zeros, construct and find the zeros of the equivalent polynomial $p(z)$, using the computed values of $s_1,\dots,s_N$ for each region.
    \item \emph{Optionally}, refine the computed root using an iterative method, such as Newton's method.
\end{enumerate}

The coefficients $\sigma_k$ of the equivalent polynomial are computed using Newton's identities, see \cite{Min2003} or Theorem 1.1.1 of~\cite{Kravanja2000}, and are given by
\begin{align}
\begin{split}
    s_1 + \sigma_1 &= 0,\\
    s_2 + s_1\sigma_1 + 2\sigma_2 &=0,\\
    & \vdots\\
    s_{N} + s_{N-1}\sigma_1 + \dots + s_1\sigma_{N-1} + N\sigma_{N}&=0,
    \label{eq:equiv_poly}
\end{split}
\end{align}
which can be solved as a linear system of equations. Once the coefficients of the equivalent polynomial $p(z)$ have been computed, all that remains is to compute its zeros. This can be done by computing the eigenvalues of the $N\times N$ companion matrix; see, for example, Lecture 25 of~\cite{Trefethen1997-ls}. In this method, repeated zeros are given by repeated eigenvalues. As mentioned in Step~4 it is possible to refine the computed zeros by using them as the initial guess for an iterative root-finding method.

The work of Delves and Lyness was built upon by Kravanja and Van Barel in~\cite{Kravanja2000,Kravanja2000ZEAL}. Specifically, they retain the idea of subdivision of the search region, but instead compute the zeros by a generalized eigenvalue problem rather than using Newton's identities and also replace the monomial basis (the use of regular moments in \cref{eq:delves}) with regular formal orthogonal polynomials (FOPs) to avoid the ill-conditioned mapping from the ordinary moments to polynomial coefficients and then to zeros. This allows a larger number of zeros in each subdivision, so the search region does not need to be subdivided as many times. Another difference is that, instead of finding repeated zeros multiple times, they are found once and their multiplicities are calculated in a separate step. For brevity, we will not detail these improvements as they are not directly relevant to the method presented in this paper and instead refer the reader to \cite{Kravanja2000,Kravanja2000ZEAL}.

While these two methods have been successfully applied to solve many problems, e.g. \cite{Cheng2025,Corless2004,Oza2019}, they require a large number of function evaluations as the quadrature method has to be applied multiple times, as shown in \Cref{sec:comp}. Furthermore, when the zeros are located close together accuracy in the computed zeros is lost, see \cite{Austin2014} for a numerical example. This motivates a different approach that retains the guarantee that the number of zeros computed in the search region is correct by using the argument principle, but instead is based on approximating the logarithmic derivative by the AAA algorithm. This improves the computation of the zeros as empirically we observe this approach to use fewer function evaluations while usually being more accurate.

\section{The AAA Algorithm}
\label{sec:AAA}

The AAA algorithm was introduced in 2018 by Nakatsukasa, S\`ete, and Trefethen in \cite{Nakatsukasa2018} and is a greedy algorithm for generating a rational approximation $r(z)$ to a real or complex-valued function $f(z)$. Despite being a relatively new algorithm, it has seen use in a wide variety of applications, many of which are summarized in \cite{AAAreview}. Typically, rational functions are more efficient than polynomials for approximating functions with singularities and discontinuities.

The AAA algorithm represents the rational approximation in a barycentric form, which, for a degree $m - 1$ rational function, is given by
\begin{equation}
    r(z) = \frac{n(z)}{d(z)} = \left.\left(\sum_{j=1}^m\frac{w_jf_j}{z-z_j}\right) \middle/ \left(\sum_{j=1}^m\frac{w_j}{z-z_j}\right)\right.,
    \label{eq:ratbary}
\end{equation}
where $z_j$ are the so-called support points, $f_j=f(z_j)$ are the function value at the support points, and $w_j$ the barycentric weights. The AAA algorithm selects the support points from a discrete set of points $Z$ according to where the error $|f-r|$ is greatest. Then, the barycentric weights are computed to minimise the least-squares error to the linearised problem $\lVert fd-n\rVert$ over $Z$. These two steps are repeated until the error $\lVert f-r \rVert$ over $Z$ is less than a specified relative tolerance. Throughout this paper we use a relative tolerance of $10^{-13}$ (the default value in \verb|Chebfun| and \verb|skzeros|), unless specified otherwise, and all computations are performed using double precision arithmetic. For more details on the steps of the AAA algorithm see \cite{Nakatsukasa2018}. Implementations of the AAA algorithm are available in MATLAB through \verb|Chebfun|~\cite{Nakatsukasa2018}, Python through the \verb|SciPy| package~\cite{Virtanen2020} (implementation written by the first author), and Julia through the \verb|RationalFunctionApproximation.jl| package~\cite{RationalFunctionApproximation.jl}.

Since the search region $\Omega$ and any subdivisions of $\Omega$ are a continuum, we use a continuum variant of the AAA algorithm described in \cite{Driscoll2024}. Here, we replace the discrete set $Z$ with the continuum $\Omega$.
Then $Z$ is adaptively chosen as a discrete subset of the boundary $\partial\Omega$ based on the existing support points. The first support point is chosen to be where the function has the greatest absolute difference from its mean value on $\partial\Omega$. Following~\cite{Driscoll2024}, $Z$ is then chosen such that there are $\max\{3, 16-m\}$ sample points between each support point along $\partial\Omega$. This ensures that there are neither too many nor too few sample points, avoiding unnecessary function evaluations or an under-resolved approximation.

Critically for our application, once the rational approximation has been computed, it is possible to find its poles and zeros. The zeros of $r(z)$ are computed by finding the zeros of $n(z)$ and can be computed by solving the following generalized eigenvalue problem~\cite{Nakatsukasa2018,klein}
\begin{equation}
    \begin{pmatrix}
        0 & w_1f_1 & w_2f_2 & \cdots & w_mf_m\\
        1 & z_1\\
        1 & & z_2\\
        \vdots & & & \ddots\\
        1 & & & & z_m
    \end{pmatrix} \mathbf{v}
    = \lambda
    \begin{pmatrix}
        0 & & & &0\\
        & 1\\
        & & 1\\
        & & & \ddots \\
        0 & & & &1
    \end{pmatrix}\mathbf{v}.
    \label{eq:eig_pol}
\end{equation}
The finite eigenvalues of \eqref{eq:eig_pol} give the zeros of $n(z)$. Likewise, the poles of the rational approximation are the zeros of $d(z)$ and can be computed by replacing $w_jf_j$ with $w_j$ in \eqref{eq:eig_pol}.

\section{Choice of Function to Approximate}
\label{sec:logderiv}

At this point, it may seem obvious that one method of solving our problem is to simply approximate $f(z)$ over the search region $\Omega$ using the AAA algorithm and then to compute the zeros of the approximation, see Section 4 of~\cite{AAAreview}. For simple problems, this approach may be sufficient. However, there are a number of improvements that can be made in order to develop a more robust approach. The first modification is to consider whether $f(z)$ is the best function to approximate. To motivate why we may wish to approximate a function other than $f(z)$ we consider a number of issues that may arise.

One issue that presents itself is the possibility for the approximation to contain spurious pole-zero pairs, known as Froissart doublets~\cite{Nakatsukasa2018}. Although we can eliminate these spurious zeros by substituting the computed zeros back into $f(z)$, it would be preferable to avoid them in the first instance.

Another issue with approximating $f(z)$ is that the computation of zeros and poles using a AAA approximation loses accuracy when the poles and zeros are of a higher order. Simple zeros can usually be computed to approximately 15 digits of accuracy. In comparison, a second-order zero can only be computed to approximately eight or nine digits of accuracy. This loss of precision can be explained by noting that the eigenvalue problem \eqref{eq:eig_pol} will have an eigenvalue of high multiplicity $\alpha$ (which can be seen to correspond to a single Jordan block of size $\alpha$), for which standard perturbation analysis shows the computed eigenvalues (when computed in a backward stable manner) will only be accurate to $\mathcal{O}\left(\epsilon^{1/\alpha}\right)$, see p.~77 of~\cite{wilkinson:1965}. Furthermore, because of this loss of accuracy, repeated zeros of $f(z)$ appear as a cluster of zeros in the approximation. This can make it difficult to distinguish between repeated zeros and clusters of zeros in $f(z)$.

To remedy these issues, we propose that, rather than using the AAA algorithm to approximate $f(z)$, it should be used to approximate the \emph{logarithmic derivative} $f^\prime(z)/f(z)$. Recall that, if $f(z)$ has a zero at $a$ of order $\alpha$ then we can write $f(z) = (z-a)^\alpha h(z)$, where the function $h(z)$ is also holomorphic. This means the logarithmic derivative is given by
\begin{equation}
    \frac{f^\prime(z)}{f(z)} = \frac{\alpha}{z-a} + \frac{h^\prime(z)}{h(z)},
\end{equation}
so the zero has been turned into a simple pole of the logarithmic derivative with residue equal to the multiplicity. The location of the poles of the rational approximation $r(z) \approx f^\prime(z)/f(z)$ can then be computed according to \eqref{eq:eig_pol} and the associated residue according to
\begin{equation}
    \res\left(r(z),z=a\right)= \left.\left(\sum_{j=1}^m\frac{w_jf_j}{a-z_j}\right) \middle/ \left(-\sum_{j=1}^m\frac{w_j}{(a-z_j)^2}\right)\right..
\end{equation}

This change to approximating the logarithmic derivative is powerful as repeated zeros become simple poles, so the loss of precision described previously is resolved. In addition, it is possible to find the multiplicity from the residue (which will be an integer $\alpha$) rather than determining which cluster of zeros makes up the repeated zero. The issue of spurious pole-zero pairs is also largely resolved, as we can be confident that the computed simple pole is correct if its residue is close to a non-zero integer, as it corresponds to the multiplicity of the zero. Note that, these benefits also apply when using AAA for pole finding.

\begin{figure}
    \centering
    \includegraphics{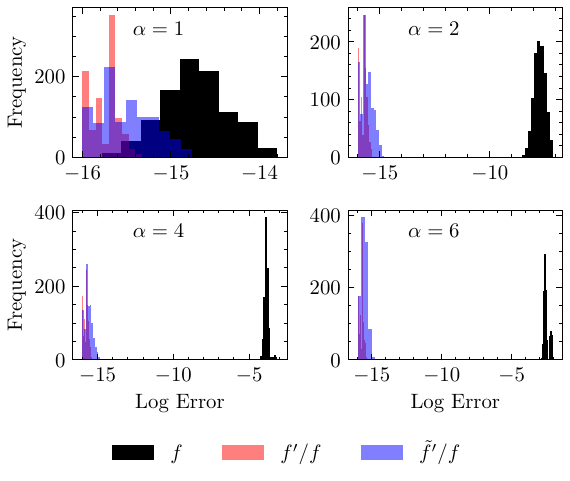}
    \caption{A comparison of the distribution of errors in the location of the zero of \eqref{eq:func_choice_test} computed by the AAA approximation to $f(z)$, $f^\prime(z)/f(z)$, and $\tilde{f}^\prime(z)/f(z)$, where $\tilde{f}^\prime(z)$ is an approximation to $f^\prime(z)$ computed using a trapezium rule approximation of Cauchy's integral formula, as given in \eqref{eq:trap_Cauchy}. Here, the zero location $a$ is a random complex number located in the unit square search region $\Omega=[0,1]+[0,1]i$. We see that the logarithmic derivative is substantially more accurate for orders $\alpha>1$.}
    \label{fig:whatfunc}
\end{figure}

The following numerical experiment was performed to verify that using the logarithmic derivative improves the accuracy of the computed zeros when using the AAA algorithm. The zero of the function
\begin{equation}
f(z)=e^z(z-a)^{\alpha},
\label{eq:func_choice_test}
\end{equation}
in the unit square search region $\Omega=[0,1]+[0,1]i=\{x+iy:x\in[0,1],y\in[0,1]\}$ was computed using AAA approximations to the original function $f(z)$ and the logarithmic derivative $f^\prime(z)/f(z)$ for 1000 different values of $a$, randomly selected from $\Omega$. This was performed for orders $\alpha=1,2,4,8$. The error in the computed zeros for the two approaches is compared in \Cref{fig:whatfunc}.
For order $\alpha=1$ using the logarithmic derivative achieves approximately one digit better accuracy although both approaches are effective. However, for order $\alpha>1$ there is a notable difference in accuracy; as predicted, using $f(z)$ to compute the zero quickly loses accuracy. In comparison, approximating the logarithmic derivative maintains around 15 digits of accuracy in the zero location. Furthermore, approximating the logarithmic derivative, we find that the number of iterations of the AAA algorithm required to reach convergence is constant as the order $\alpha$ increases. In contrast, when approximating the original function, the number of iterations grows with the order, increasing the runtime of the AAA algorithm. This is another benefit of approximating the logarithmic derivative instead of the original function.

One further question to consider is: what if $f^\prime(z)$ is not explicitly known? One obvious idea might be to use a different formulation of the logarithmic derivative. So far we have only considered $f^\prime(z)/f(z)$ yet if $f(z)$ is holomorphic and non-zero then we have the following result
\begin{equation}
    \frac{\diff{}}{\diff{z}}\log\left(f(z)\right)=\frac{f^\prime(z)}{f(z)}.
\end{equation}
This suggests an alternative approach where we first compute an AAA approximation to $\log(f(z))$, then differentiate the AAA approximation using the formula given in Proposition~11 of~\cite{Schneider1986}, and finally compute its poles. However, the reason why this alternative formulation is not used is that taking $\log(f(z))$ introduces branch cuts to the function we have to approximate. This is problematic as the discontinuity at the branch cut is difficult for the AAA algorithm to approximate. This is discussed in more detail in \Cref{sec:brach}. Since different branches of the complex logarithm differ by a constant, they have the same derivative. Hence, branch cuts are not an issue in the $f^\prime(z)/f(z)$ formulation of the logarithmic derivative, assuming that $f(z)$ itself is free from branch cuts. See \Cref{sec:brach} for handling cases where this is not true. Another possible derivative-free approach is to use two AAA approximations: first we approximate $f(z)$ by $r(z)$ and then we approximate $r^\prime(z)/f(z)$, where $r^\prime(z)$ is the derivative of the inner approximation is computed using the formula given in. This approach is discussed in Section 6 of~\cite{AAAreview}. In practice, we find that this formula is too numerically unstable near the support points of the inner approximation for this to be a reliable approach.

In our experience, the most robust derivative-free approach is based on a discretisation of Cauchy's integral formula, as discussed in~\cite{Austin2014,Bornemann2010,Lyness1967b,Lyness1967}. Recall that, if $\gamma$ is a simple closed curve containing $z$ and $f(z)$ is holomorphic in the interior of $\gamma$, then Cauchy's integral formula tells us that the derivative of $f(z)$ can be computed as
\begin{equation}
    f^\prime(z) = \frac{1}{2\pi i}\oint_{\gamma} \frac{f(\xi)}{(\xi-z)^2} \diff \xi.
\end{equation}
Making a change of variables so the pole is centred at the origin and taking $\gamma$ to be a circular contour of radius $r$ about it yields
\begin{equation}
    f^\prime(z) = \frac{1}{2\pi r} \int_0^{2\pi} e^{-i\theta}f\left(z + re^{i\theta}\right)\diff \theta.
\end{equation}
We can then exploit the exponential convergence of the trapezium rule for holomorphic functions over periodic intervals~\cite{davis1959,Trefethen2014} and approximate $f^\prime(z)$ by
\begin{equation}
    \tilde{f}^\prime(z) = \frac{1}{mr}\sum_{j=0}^{m-1} e^{-2\pi ij/m}f\left(z + re^{2\pi ij/m}\right),
    \label{eq:trap_Cauchy}
\end{equation}
where $m$ is the number of nodes used. Specifically, we use the algorithm given in Figure 3 of~\cite{Bornemann2010} in which $m$ is selected by a successive doubling procedure.

Now we have a means of accurately approximating the derivative, we can compute the logarithmic derivative without explicitly knowing $f^\prime(z)$. \Cref{fig:whatfunc} shows the accuracy of the computed zeros of \eqref{eq:func_choice_test} when computed with a AAA approximation of $\tilde{f}^\prime(z)/f(z)$, where we use a radius of $10^{-2}$ and a relative tolerance of $10^{-15}$ when selecting $m$. Furthermore, as we are computing a AAA approximation of an approximation we loosen the relative tolerance of the continuum AAA algorithm to $10^{-12}$. We observe that very little accuracy is lost compared to when the true derivative is known. Therefore we believe this to be a viable approach for when the derivative $f^\prime(z)$ is not readily available, though at the expense of extra function evaluations.

\section{Subdivision}
\label{sec:subdiv}

As demonstrated in \Cref{fig:whatfunc}, using the AAA algorithm to approximate the logarithmic derivative can be an accurate method to compute the zeros in a region. However, the number of iterations of the AAA algorithm is equal to the degree of the rational approximation, which in turn must be greater than or equal to the number of zeros (poles of the logarithmic derivative) in the region if it is to locate them all. This means if there are many zeros in the region then many iterations of the AAA algorithm will be necessary to approximate the logarithmic derivative. This is problematic as the runtime of the continuum AAA algorithm is quartic in the number of iterations performed. Additionally if $f(z)$ has too many zeros, the AAA algorithm may not be able to converge to the requested tolerance. This motivates the following divide-and-conquer type approach: when there are too many zeros in a region we will subdivide the region to make it easier to approximate with the AAA algorithm.

\begin{figure}
    \centering
    \includegraphics{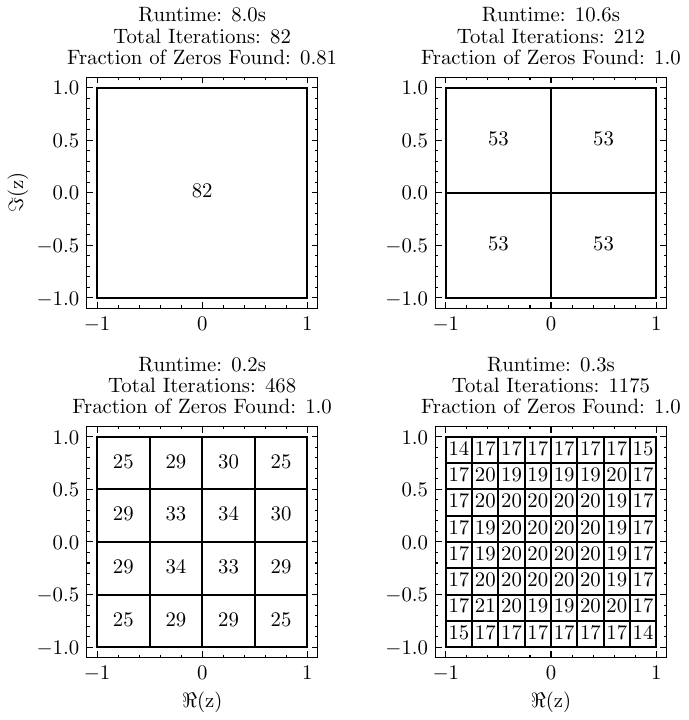}
    \caption{Numerical experiment showing the effect of subdivision on the runtime and the total number of iterations of the AAA algorithm as well as the proportion of the zeros of \eqref{eq:test_func} found as the search region $\Omega=[-1,1]+[-1,1]i$ is subdivided. The number within each square denotes the number of iterations required to approximate the logarithmic derivative in that region. Despite increasing the total number of iterations required, subdivision decreases the total runtime. The timings were computed on an unladen laptop with 16 gigabytes of memory and an 11th Gen Intel\textregistered \ Core\texttrademark \ i7-11800H processor.}
    \label{fig:AAA_subdivision}
\end{figure}

To test this idea, we perform a numerical experiment inspired by \cite{sub}. Given the grid of points $A=\left\{(x+iy)/9: (x,y)\in\{-9,-7,-5,-3,-1,1,3,5,7,9\}^2\right\}$, we consider the function
\begin{equation}
    f(z) = \prod_{a\in A}(z-a),
    \label{eq:test_func}
\end{equation}
i.e. $f(z)$ has 100 zeros. We wish to compute the zeros of \eqref{eq:test_func} by approximating its logarithmic derivative
\begin{equation}
    \frac{f^\prime(z)}{f(z)} = \sum_{a\in A} \frac{1}{z-a},
    \label{eq:test_func2}
\end{equation}
with the AAA algorithm. We begin by approximating \eqref{eq:test_func2} with the AAA algorithm on the region $\Omega=[-1,1]+[-1,1]i$, recording the number of iterations of the AAA algorithm required, the runtime taken to generate the approximation, and the proportion of the zeros found. Then, we successively divide the region into quarters, recording the number of iterations of the AAA algorithm required on each subregion and the total runtime. The results of this experiment are shown in \Cref{fig:AAA_subdivision}. Here, we observe that, while subdividing the region leads to a higher number of total iterations of the AAA algorithm, the runtime can be greatly decreased. Furthermore, we see that initially not all the zeros of \eqref{eq:test_func} are found but after subdivision they are all computed. However, we also see that subdividing the region too many times leads to the runtime increasing again, as clearly more problems need to be solved. This experiment motivates the method described in \Cref{sec:AAA_Quad}, which combines AAA approximation of the logarithmic derivative and subdivision to find the zeros of a function even when there are many zeros. To our knowledge, the idea of combining subdivision with the AAA algorithm has only been considered once before in the literature in \cite{AAAsubdivision} in a different context (where subdivision is observed to improve the accuracy).

\section{Functions with Branch Cuts}
\label{sec:brach}

\begin{figure}
    \centering
    \includegraphics{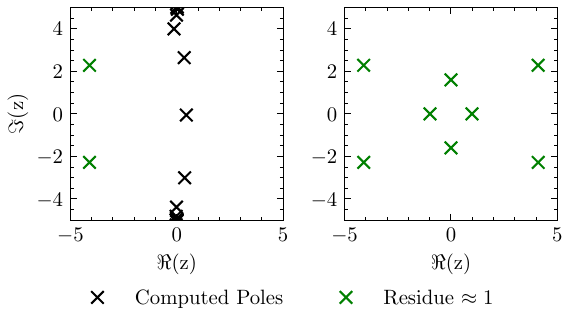}
    \caption{ Computed poles of the AAA approximation of the logarithmic derivative of \eqref{eq:branch} (left) and the logarithmic derivative of the product of the Riemann sheets \eqref{eq:sheet1} and \eqref{eq:sheet2} (right) in the region $\Omega=[-5,5]+[-5,5]i$. (Left Plot) Here we see clusters of poles at the intersection of the boundary of $\Omega$ and the branch cuts. Here, AAA terminates without achieving the desired relative tolerance for the approximation. (Right Plot) By taking the product of Riemann sheets the function is free from branch cuts. Thus, AAA successfully converges and we are able to identify all the zeros in $\Omega$.}
    \label{fig:branch_cut}
\end{figure}

One remaining issue to deal with is the possibility that the function may have branch cuts in or near the search region $\Omega$, which will have to be approximated by the AAA algorithm. For example, consider the function
\begin{equation}
    f(z) = \sin\left(\sqrt{z^2 + 1}\right) - z,
    \label{eq:branch}
\end{equation}
from \cite{Kowalczyk2017}. For $f(z)$ to be single-valued and assuming the standard convention that the square root has a branch cut along the negative real axis, $f(z)$ must have branch cuts along $(i,\infty)$ and $(-i,-\infty)$. The poles of the AAA approximation of the logarithmic derivative of $f(z)$ in the search region $\Omega=[-5,5]+[-5,5]i$ are shown in \cref{fig:branch_cut}. Here we encounter two problems. First, the AAA algorithm approximates the discontinuity at the branch cuts with a cluster of interlaced zeros and poles at the boundary of $\Omega$, see Section~7.4 of~\cite{Trefethen2023}. Although we can still identify some of the zeros of $f(z)$ by looking for the poles with residue approximately equal to 1, the runtime of the AAA algorithm is greatly increased as many iterations are required to insert enough poles to approximate the discontinuity at the branch cuts. Second, the zero at $z\approx1.5999i$ has been missed.

In order to develop an approach that works around these limitations, we use the result from Section~2 of~\cite{Kowalczyk2017} and take the pointwise product of the Riemann sheets of $f(z)$. We have two Riemann sheets
\begin{align}
    f^+(z) &= \sin\left(\sqrt{z^2 + 1}\right) - z,
    \label{eq:sheet1}\\
    f^-(z) &= -\sin\left(\sqrt{z^2 + 1}\right) - z,
    \label{eq:sheet2}
\end{align}
so we form the pointwise product $F(z)=f^+(z)\cdot f^-(z)$. Again, we assume the standard convention for the branch cut of the square root. The function $F(z)$ is now much easier to compute a AAA approximation of because it no longer has any discontinuities. Plotting the poles of the logarithmic derivative of $F(z)$, see \Cref{fig:branch_cut}, we see that all the zeros in the region have been found successfully. Then it is simply a task of determining which Riemann sheet they belong to by evaluating $f^+(z)$ and $f^-(z)$ at the computed zeros to see which function is equal to zero. This approach allows us to compute the zeros of functions using the AAA approximation of the logarithmic derivative even when $f(z)$ has branch cuts.

\section{Comparison of AAA with Quadrature-Based Methods}
\label{sec:comp}

In this section, we compare the accuracy of computing the zeros of a AAA approximation of the logarithmic derivative against the prior methods presented in \Cref{sec:quad_methods}. These are the methods of Delves and Lyness, Kravanja and Van Barel with a monomial basis, and Kravanja and Van Barel with a FOP basis. The implementation of the former two methods is provided in the repository linked in \Cref{sec:code} and the implementation of the later method is provided by the \verb|cxroots| package~\cite{cxroots}. The performance of the methods is compared by looking at the error in the computed zeros against the number of function evaluations of the logarithmic derivative required. For the quadrature-based methods the number of function evaluations is controlled by changing the relative tolerance of the quadrature. For the AAA algorithms we change the relative tolerance that terminates the algorithm. To test the methods, we consider finding an increasing number of zeros in the unit square search region $\Omega=[0,1]+[0,1]i$. We consider all methods without any additional subdivision as all four methods can be used with the same subdivision strategy. Instead, this experiment is designed to test the performance after the subdivision has already been completed. Additionally, we consider the performance of the methods without any form of polishing with an iterative method. This is because, with polishing, all methods can usually achieve approximately 16 digits of accuracy through the use of additional function evaluations. Instead, we wish to understand the performance of the methods before any polishing steps.

\begin{figure}
    \centering
    \includegraphics{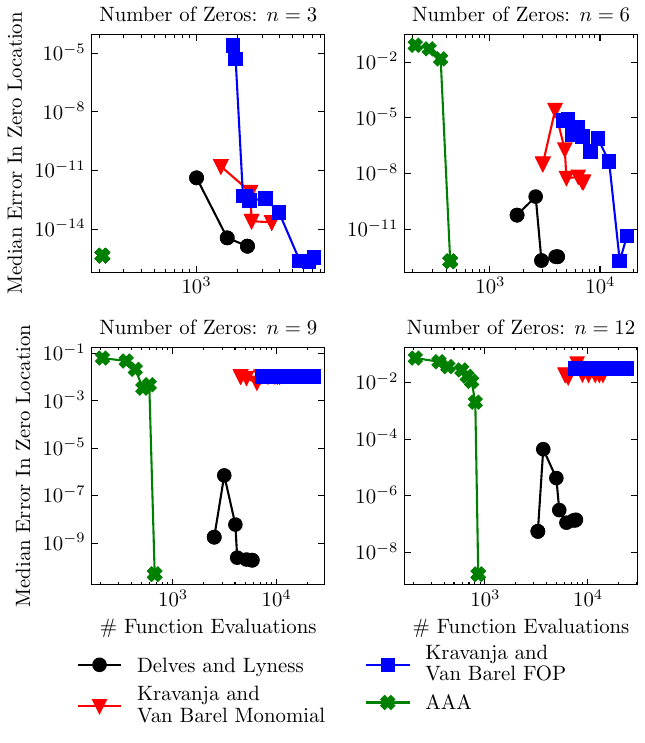}
    \caption{Error in the computed zeros of \eqref{eq:compfunc} as the number of zeros $n$ is increased. Note, that for the $n=9,12$ cases the two Kravanja and Barel methods fail to give a meaningful solution.}
    \label{fig:comparison_number_zeros}
\end{figure}

To compare the methods, we consider the accuracy of the computed zeros against the number of function evaluations required to compute them for the test function
\begin{equation}
    f(z;n)=\prod_{j=0}^n\left(z-a_j(n)\right),
    \label{eq:compfunc}
\end{equation}
as the number of zeros $n$ increases, where $a_j(n)=(0.1 + 0.8j/n) + 0.5i$, i.e. there is a line of $n + 1$ zeros across the middle of the square. The results of this experiment are shown in \Cref{fig:comparison_number_zeros}. Here, the AAA algorithm significantly outperforms the other methods, using far fewer function evaluations while achieving equal or better accuracy. For the $n=9,12$ cases we see that the two Kravanja and Van Barel methods fail to converge. In their full algorithm, this would mean that this region would have to be subdivided further. This numerical experiment shows that the AAA approximation of the logarithmic derivative is a good candidate to replace the quadrature-based approximation techniques used to find the zeros of the subregions as it can achieve equally good accuracy while using far fewer function evaluations.

\section{Combining AAA and Quadrature-Based Subdivision}
\label{sec:AAA_Quad}

In this section, we introduce a new method based on subdivision and AAA approximation of the logarithmic derivative, motivated by the earlier numerical experiments. Given a user-supplied holomorphic function $f(z)$, derivative $f^\prime(z)$, and rectangular search region $\Omega\subset\mathbb{C}$, the method can be summarized at a high level by the following two steps:
\begin{enumerate}
    \item \emph{Subdivision}: compute the number of zeros $N$ in the search region $\Omega$ using Cauchy's argument principle. If $N$ is greater than a user-defined threshold $M$ then there are too many zeros in the region, so the region is subdivided and this step is repeated with the two child regions. Here $M$ must be greater than the highest multiplicity of any zero, otherwise the subdivision will never complete unless a limit on the number of iterations is set.
    \item \emph{Approximation}: Once the number of zeros $N$ in each region is less than or equal to the user defined threshold $M$, we approximate the logarithmic derivative $f^\prime(z)/f(z)$ in each subregion using the continuum AAA algorithm. We then compute the poles of the approximation as these correspond to the zeros of the original function if the associated residue is close to a positive integer.
\end{enumerate}
We now proceed to explain each of these steps in more detail.

\subsection{Subdivision}

\begin{figure}
    \centering
    \includegraphics{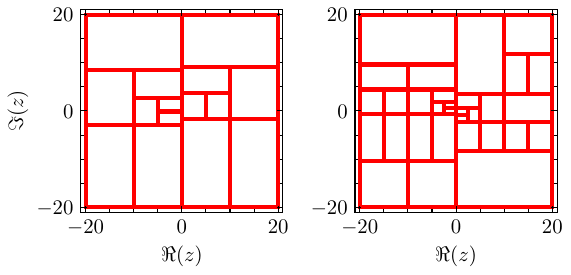}
    \caption{Subdivision of the function $f(z)=z^{50}+z^{12}-5\sin(20z)\cos(12z) - 1$ within the search region $\Omega=[-20,20] + [-20,20]i$. In the left plot the maximum number of zeros in any given region is $M=50$ and in the right $M=25$. The asymmetry of the subdivisions is due to the fact that when a zero falls on or close to an edge we perturb the edge by a small random amount to avoid repeatedly placing an edge through a zero.}
    \label{fig:subdiv_examples}
\end{figure}

The search region $\Omega\subset\mathbb{C}$ is the area of the complex plane in which the method looks for zeros of $f(z)$. For ease of exposition and implementation we will assume $\Omega$ to be rectangular, however, it would be possible to generalise this method to the subdivision of other shaped regions. For example, \cite{Delves1967} also includes an approach for subdividing a disk region. The subdivision stage begins by computing the total number (counted according to their multiplicity) of zeros in $\Omega$ by Cauchy's argument principle
\begin{equation}
    N\left(f, f^\prime, \partial\Omega\right) = \frac{1}{2\pi i}\oint_{\partial\Omega} \frac{f^\prime(z)}{f(z)}\diff{z}.
    \label{eq:num_zeros}
\end{equation}
The value of the integral is computed with numerical quadrature. We use a scheme based on 10-point Gauss and 21-point Kronrod rules with globally adaptive interval subdivision~\cite{deDoncker1978}.

The method then proceeds based on the computed result of Cauchy's argument principle around each rectangle. There are two possibilities that need to be addressed. The first is that the quadrature completes successfully and there are too many zeros in the region, i.e. $N>M$. This means the region must be subdivided. To do this, we divide the rectangular region into two regions of equal area by subdividing either vertically or horizontally. We choose the subdivision that minimises the aspect ratio of the two resulting rectangular regions. We denote subregions that are created together by splitting a region, siblings, and the region that was split, the parent region. The subdivision process repeats recursively on each of the two subregions: the number of zeros in each of the two subregions is computed and then subdivision is performed again and again, as required. Note that, when a region is divided in half the resulting quadrature along the shared edge is equal up to a difference in sign so this result can be reused.

The second possibility is that there is a zero close to or on the boundary $\partial\Omega$, in which case the numerical quadrature will fail as \eqref{eq:num_zeros} is divergent or slow to converge. We are able to detect this as the quadrature routine warns that it is unable to achieve the desired integration tolerance. This approach is also used in \cite{Kravanja2000ZEAL} to detect when there are poles close to or on the integration contour. If this occurs on the first iteration we raise an error to notify the user that there is a zero on the boundary of the search region as this cannot be resolved without changing the supplied search region $\Omega$. If this occurs in an iteration after the region has been subdivided then it must be the inserted edge that is faulty as it is already known that the boundary is free from zeros. In this case, the edge dividing the region is perturbed and the quadrature attempted again. Care must be taken to also adjust any sibling regions that share the perturbed edge, after which the two regions are no longer of equal area.

Finally, once $N\leq M$ is satisfied in every subregion, the subdivision algorithm terminates at which point we have that
\begin{equation}
    \Omega = \bigcup\limits_{\Omega_i\in A} \Omega_i,
\end{equation}
where $A$ is the set of accepted regions. The steps described here are also summarized by the pseudocode given in \Cref{alg:subdiv}. \Cref{fig:subdiv_examples} shows examples of subdivided search regions.

\begin{algorithm}
    \caption{Subdivision.}
    \label{alg:subdiv}
    \textbf{Input:} A holomorphic function $f(z)$ and its derivative $f^\prime(z)$; a search region $\Omega$; and $M$ the maximum value of the argument principle in each subregion.

    \textbf{Output:} A subdivision of $\Omega$ such that the value of the argument principle is less than or equal to $M$ in each subregion.
    \begin{algorithmic}[1]
        \STATE Initialise the queue of regions that may require subdivision $Q = \{\Omega\}$ and the set of accepted regions $A=\emptyset$.
        \WHILE{$Q \neq \emptyset$}
            \STATE Retrieve the next region $\Omega_i$ from the front of the queue $Q$.
            \STATE Compute the number of zeros $N(f,f^\prime,\Omega_i)$ using the argument principle.
            \IF{Quadrature fails to converge}
                \IF{$\Omega_i = \Omega$}
                    \STATE Output an error as there is a zero on or close to $\partial\Omega$ (the user supplied region).
                \ENDIF
                \STATE Remove any siblings of $\Omega_i$ that share the same edge for which the quadrature failed to converge from $Q$.
                \STATE Reattempt to subdivide the parent region of $\Omega_c$ by adding a perturbation to the failed edge and add the new regions to the queue $Q$.
            \ELSIF{$N$ is approximately a positive integer}
                \IF{$N>M$}
                    \STATE Subdivide $\Omega_i$ and add the resulting regions to the queue $Q$.
                \ELSE
                    \STATE Add $\Omega_i$ to the set of accepted regions $A$.
                \ENDIF
            \ELSE
                \STATE Output an error as a non-integer argument principle may mean that the function is not holomorphic in $\Omega$ or that the quadrature has failed to converge.
            \ENDIF
        \ENDWHILE
    \end{algorithmic}
\end{algorithm}

\subsection{Approximation}

Once the subdivision of the search region is complete, the method proceeds to perform AAA approximation on each of the subregions where Cauchy's argument principle tells us zeros are present, i.e. $N$ is greater than zero. We apply the continuum AAA method described in \Cref{sec:AAA}, by parametrising and then sampling the boundary of the subregion to approximate the logarithmic derivative $f^\prime(z)/f(z)$. The poles and associated residues are then computed using the generalised eigenvalue problem \eqref{eq:eig_pol}. Then, we discard any poles that lie outside of the subregion under the assumption that we will locate them more accurately when we apply AAA to the subregion that contains them. Next, we check the residue of the remaining poles. We know that the residue gives the multiplicity of the zero, hence, it should take (up to a tolerance of say $10^{-2}$) a positive integer value. Consequently, we discard any poles where this is not the case as they are likely spurious.

Once the zeros of a region have been computed, the method compares the sum of the multiplicities of the zeros against the expected value given by Cauchy's argument principle. If they match then we can be confident that we have correctly found the right number of zeros. If they do not match, we assume that the AAA approximation is inaccurate and we subdivide the region again in order to obtain a better approximation. This generally only occurs if $M$ is set too high leading to insufficient subdivision of the search region. We repeat this until the result computed from the AAA approximation and Cauchy's argument principle are in agreement.

The method concludes when this step has successfully been applied to all subregions, at which point we have found all the zeros of $f(z)$ in $\Omega$. The process described here is ``embarrassingly parallel'' as we can treat each subregion independently. If required, the user may refine the accuracy of the computed zeros by using them as a starting guess for an iterative method, such Newton's or Muller's method. This section is summarized by the pseudocode given in \Cref{alg:Zerofinding}.

\begin{algorithm}
    \caption{Zero finding with AAA and subdivision.}
    \label{alg:Zerofinding}
    \textbf{Input:} A holomorphic function $f(z)$ and its derivative $f^\prime(z)$; a search region $\Omega$; and $M$ the maximum value of the argument principle in each subregion.

    \textbf{Output:} The zeros of $f(z)$ that lie in $\Omega$ and their associated multiplicities.
    \begin{algorithmic}[1]
        \STATE Initialise the queue of accepted regions $A$ as the output of \Cref{alg:subdiv}, and the set of located zeros and their associated multiplicities $Z = \emptyset$.
        \WHILE{$A\neq\emptyset$}
            \STATE Retrieve the next region $\Omega_i$ from the front of the queue $A$.
            \IF{number of zeros $N(f,f^\prime,\Omega_i) > 0$}
                \STATE Use the continuum AAA algorithm to generate the rational approximation $r(z)\approx f^\prime(z)/f(z)$ on $\Omega_i$.
                \STATE Compute the poles $a_i$ and associated residues $\alpha_i$ of $r(z)$.
                \STATE Discard any poles $a_i$ which lie outside of $\Omega_i$.
                \STATE Discard any poles $a_i$ for which the residue $\alpha_i$ that is not approximately a positive integer as they are spurious.
                \IF{$\sum_i\alpha_i \approx N(f,f^\prime,\Omega_i)$}
                    \STATE Add each $(a_i, \alpha_i)$ to $Z$.
                \ELSE
                    \STATE Subdivide $\Omega_i$.
                    \IF{quadrature for computing the number of zeros using the argument principle of each of the children of $\Omega_i$ converges}
                        \STATE Add the children of $\Omega_i$ to the queue $A$.
                    \ELSE
                        \STATE Resubdivide $\Omega_i$ adding a perturbation to the subdivision as there is a zero near or on a boundary.
                        \STATE Go to 14.
                    \ENDIF
                \ENDIF
            \ENDIF
        \ENDWHILE

    \end{algorithmic}
\end{algorithm}

\section{Numerical Experiments}
\label{sec:Num_expr}

In this section we demonstrate the method on three test problems to verify it works as intended.

\subsection{Example 1: Polynomial with Quasi-Random Zeros in the Unit Square}

\begin{figure}
    \centering
    \includegraphics{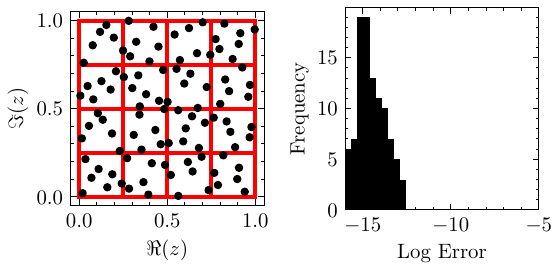}
    \caption{Results for \eqref{eq:test_func_1} with $M=7$. (Left Plot) The location of the computed zeros and subdivision of the search region $\Omega$. (Right Plot) The distribution of the error in computed zeros.}
    \label{fig:quasi-random}
\end{figure}

The first example we consider is $N$ quasi-randomly chosen simple zeros in the search region $\Omega=[0,1]+[0,1]i$. That is
\begin{equation}
    f(z) = \prod_{j=1}^N(z-z_j),
    \label{eq:test_func_1}
\end{equation}
where $z_j$ is quasi-randomly chosen in $\Omega$ using a Sobol' sequence~\cite{Sobol1967}. We generate the points quasi-randomly as it avoids clusters of arbitrarily close zeros which may be infeasible to subdivide. While this function is not motivated by a practical application, this allows us to quantify the accuracy of the computed zeros.

Using the method presented in \Cref{sec:AAA_Quad} we can compute all 100 zeros of \eqref{eq:test_func_1} in approximately 1 second on our laptop. We initially subdivide the region until there are at most $M=7$ zeros in each subregion. \Cref{fig:quasi-random} shows how the domain was subdivided by the method and the distribution of the errors in the computed zeros. Here, we see that all the zeros have been successfully computed to at least 12 decimal places of accuracy.

\subsection{Example 2: Transcendental Equation Modelling the Stability of an Annular Combustion Chamber}
\label{sec:eigenvalues}

\begin{figure}
    \centering
    \includegraphics{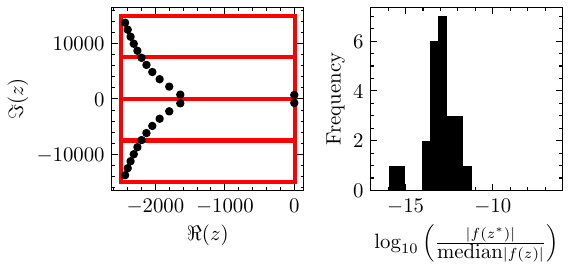}
    \caption{Results for \eqref{eq:annular} with $A = -0.19435$, $B = 1000.41$, $C = 522463$, $T = 0.005$, and $M=7$. (Left Plot) The location of the computed zeros and subdivision of the search region. (Right Plot) The distribution of the relative residuals of the computed zeros. Here $z^*$ denotes the computed zero and $\mbox{median}|f(z)|$ is the median value of $|f(z)|$ in $\Omega$.}
    \label{fig:example2}
\end{figure}

The second example we consider is from Section 4.2 of~\cite{Dellnitz2002} and relates to a problem from the Corporate Technology Department of Siemens (Munich) in which they model the stability of a flow inside an annular combustion chamber. Analysing this model leads to the following zero-finding problem
\begin{equation}
    f(z) = z^2 + Az + B\exp(-Tz) + C,
    \label{eq:annular}
\end{equation}
where $A, B, C, T \in\mathbb{R}$ are constants. We can compute the zeros $\eqref{eq:annular}$ in the region $\Omega=[-2500, 10] + [-15000,15000]i$ in less than a second, the results of which are shown in \Cref{fig:example2}. As zeros of \eqref{eq:annular} are not known, we cannot check the error in the computed zeros so instead we check the magnitude of relative residual.

\subsection{Example 3: Eigenvalues of a Matrix Computed Two Ways}
For the final example, inspired by Section 5 of~\cite{Austin2014}, we consider how we can use the method presented in this paper to compute the eigenvalues of a matrix using two different approaches. Specifically, we will consider the $50\times 50$ circulant matrix $A$ with real entries uniformly sampled from the set $\{-0.4, 0.4\}$.

\subsubsection{Zero finding}

\begin{figure*}
    \centering
    \begin{subfigure}[t]{\textwidth}
        \centering
        \includegraphics{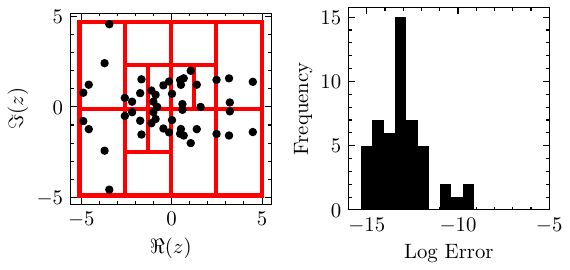}
        \caption{$f(z)=\det\left(A-zI\right)$}
        \label{fig:example3a}
    \end{subfigure}
    \begin{subfigure}[t]{\textwidth}
        \centering
        \includegraphics{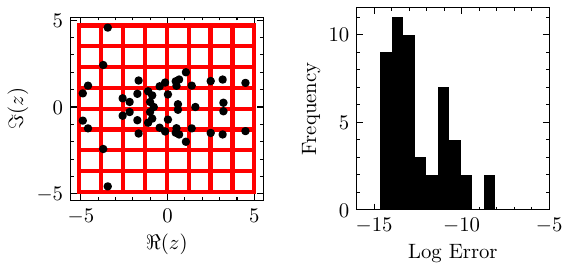}
        \caption{$f(z) = u^*\left(A- zI\right)^{-1}v$}
        \label{fig:example3b}
    \end{subfigure}
    \cprotect\caption{(Left Plots) The location of the computed eigenvalues and subdivision of the search region. (Right Plots) The distribution in the error of the computed eigenvalues compared against those computed by the \verb|dgeev| LAPACK routine.}
\end{figure*}

The first way we will approach computing the eigenvalues of $A$ is as a zero-finding problem, as discussed in Section 5 of~\cite{Austin2014}. From elementary linear algebra we know that the eigenvalues of $A$ are the zeros of the function
\begin{equation}
    f(z) = \det\left(A - zI\right),
    \label{eq:3a}
\end{equation}
for which the derivative can be computed using Jacobi's formula and is given by
\begin{equation}
    f^\prime(z) = -\tr\left(\adj \left(A-zI\right) \right),
\end{equation}
where $\tr(\cdot)$ denotes the trace and $\adj(\cdot)$ denotes the adjugate matrix. \Cref{fig:example3a} shows the results for \eqref{eq:3a} with the search region $\Omega=[-5.1, 5] + [-4.9,4.7]i$. As per the previous examples, we set $M=7$. Here $\Omega$ is deliberately asymmetric about the real and imaginary axes as otherwise subdivision would occur along the real axis where some eigenvalues are located and the quadrature would fail. As previously discussed the method can handle this automatically by perturbing the subdivision, however, here we deliberately prevent this from occurring to save computation. The error distribution plot in \Cref{fig:example3a} shows that we are successfully able to compute all of the eigenvalues of $A$ to at least 10 decimal places, often 14.

It is worth noting, however that this approach is difficult to use for large-scale problems, as the determinant is not easy to compute with high relative accuracy and is prone to underflow and overflow.

\subsubsection{Pole finding}

The second way of computing the eigenvalues of $A$ we will consider is based on \cite{Austin2014,Austin2015,Sakurai2003} and uses the scalarised resolvent
\begin{equation}
    f(z) = u^*\left(A- zI\right)^{-1}v,
    \label{eq:scalarised_resolvent}
\end{equation}
where $u$ and $v$ are random vectors. With probability 1 \eqref{eq:scalarised_resolvent} has poles at the eigenvalues of $A$ so, unlike the other examples we have considered, this is a pole-finding problem. The derivative of the scalarised resolvent $f(z)$ is given by
\begin{equation}
    f^\prime(z) = -u^*\left(\left(A- zI\right)^{-1}\right)^2v.
\end{equation}
Note that, here the function $f(z)$ is meromorphic \emph{not} holomorphic and it is the location of the poles rather than the location of the zeros that are the quantity of interest. This poses an added difficulty as the argument principle \eqref{eq:num_zeros} is now equal $N - P$, the difference between the number of zeros and poles. This is problematic as this means the number of zeros and poles cancel each other out so the automatic subdivision of the search region and the verification that the correct number of poles are computed both fail. Instead, in order to apply the ideas of this paper we simply manually specify how many times the search region should be subdivided and then apply AAA approximation of the logarithmic derivative as normal. The results for $\Omega=[-5.1, 5] + [-4.9,4.7]i$ with 6 levels of subdivision are shown in \Cref{fig:example3b}. This example stands to show how the ideas from this paper could be adapted to finding the poles and zeros of a meromorphic function.

\section{Summary and Discussion}
\label{sec:conclusion}
We have introduced a method for computing all the zeros of a holomorphic function in a specified search region by combining the flexibility of the AAA algorithm with the guarantees ensured by use of the argument principle. By approximating the logarithmic derivative rather than the original function we have been able to resolve the issues presented by higher order zeros and spurious pole-zero pairs. We have also shown how a divide-and-conquer approach using subdivision is effective at dramatically reducing the runtime required by AAA. We hope that the implementation provided  will be used to solve a range of problems. To conclude we state some limitations of the method and avenues for further research.

As stated in the introduction, we have assumed $f(z)$ to be holomorphic in $\Omega$. This is because if $f(z)$ were instead meromorphic in $\Omega$, as is the case for the example function \eqref{eq:scalarised_resolvent}, then the argument principle gives
\begin{equation}
    \frac{1}{2\pi i}\oint_{\partial\Omega} \frac{f^\prime(z)}{f(z)}\diff{z} = N - P,
\end{equation}
where $P$ is the number of poles and, as before, $N$ is the number of zeros in $\Omega$, both counted according to the multiplicity. Here, the issue that arises is that if there are both poles and zeros in the same region then we can no longer uniquely determine how many poles and how many zeros there are. This poses an issue as subdivision may halt prematurely and we can no longer have a reference to check the number of poles and zeros computed by AAA against. Of course, if $f(z)$ is meromorphic and non-zero in $\Omega$ then this is not an issue and the method can successfully compute all of the poles of $f(z)$. The only difference is that the residue of the simple poles of the logarithmic derivative computed by AAA should now be a negative integer rather than a positive integer. However, if $f(z)$ is meromorphic and we do not have this extra restriction, while there is nothing to stop the method finding some of the poles and zeros of $f(z)$, we can no longer guarantee they will all be found. In practice, the example in \Cref{sec:eigenvalues} shows that we can still successfully compute all the poles in this case using manual subdivision. Further work is required to make the automatic subdivision and verification that the correct number of poles and zeros has been computed robust in the more general case of $f(z)$ being meromorphic.

The use of the continuum AAA algorithm closely follows that of the method described in \cite{Driscoll2024}, only sampling the function on the boundary of the region. In our experiments this has proved effective at generating an accurate rational approximation without using excessive function evaluations. However, we have observed experimentally that poles closer to the boundary, and thus where the function is evaluated, are computed to a higher accuracy than those in the center of the region. This motivates further research to develop a variant of the continuum AAA algorithm that adaptively samples the function not only on the boundary but also within the region. Furthermore, other variants of the AAA
algorithm could be investigated, such as those introduced in~\cite{AAAvarient}.

The AAA algorithm is a relatively new method and so far not much has been proved about its convergence. We hope that as these results are developed they may provide some insight into when the method may possibly fail. In the meantime, the use of the argument principle at least ensures that we are able to detect such failures, if they occur.

With progress in these three avenues we hope to further improve the accuracy, performance, and robustness of the method presented.

\section*{Acknowledgments}
We thank Ellen O'Carroll, Sam Harris, and Nick Trefethen for their helpful comments. JB was funded by an MMSC Industrially Funded Scholarship from the Mathematical Institute, University of Oxford and by an EPSRC DTP Scholarship, provided by the University of Bristol.

\bibliographystyle{siamplain}
\bibliography{references}
\end{document}